\documentclass[11pt,twoside]{amsart}
\usepackage{mathrsfs}
\usepackage{amssymb}
\usepackage{amsmath}
\usepackage[all]{xy}
\usepackage{amsthm}
\numberwithin{equation}{section}

\newtheorem{theorem}{Theorem}[section]
\newtheorem{proposition}[theorem]{Proposition}
\newtheorem{lemma}[theorem]{Lemma}

\newtheorem{corollary}[theorem]{Corollary}
\newtheorem{question}[theorem]{Question}
\newtheorem{theorem*}{Theorem}

\theoremstyle{definition}
\newtheorem{definition}[theorem]{Definition}

\theoremstyle{remark}
\newtheorem{remark}[theorem]{Remark}

\usepackage{paralist}

\title[Remarks on Gorenstein weak injective and weak flat modules]{Remarks on Gorenstein weak injective \\ and weak flat modules}


\author[T. Zhao, Y. Xu]{Tiwei Zhao, Yunge Xu}
\address[T. Zhao]{School of Mathematical Sciences, Qufu Normal University,  273165 Qufu,  P. R. China}
\email{tiweizhao@qfnu.edu.cn}
\address[Y. Xu]{Faculty of Mathematics and Statistics, Hubei University, 430062 Wuhan, P. R.  China}
\email{xuy@hubu.edu.cn}

\begin{document}

\baselineskip=16pt


\maketitle

\begin{abstract}
In this paper, we introduce the notions of Gorenstein weak injective and weak flat modules respectively in terms of weak injective and weak flat modules, which is larger than  classical classes of Gorenstein injective and flat modules. In this new setting, we  characterize rings over which all modules are Gorenstein weak injective. Moreover, we also  discuss a relation between weak cosyzygy and Gorenstein weak cosyzygy of a module, and the stability of Gorenstein weak injective modules.
\end{abstract}

\section{Introduction}

Throughout $R$ is an associative ring with identity and all modules are unitary. Unless stated otherwise, an $R$-module
will be understood to be a left $R$-module. Given an $R$-module $M$, we denote by $pd_R(M)$ and $fd_R(M)$
the projective and flat dimensions respectively, and by the character module  $M^+=\mbox{Hom}_\mathbb{Z}(M,\mathbb{Q}/\mathbb{Z})$. For unexplained concepts and notations, we refer the
readers to \cite{EJ,Ro}.

In 1970, Stenstr\"{o}m  \cite{St} introduced the notion of FP-injective modules, and generalized the homological properties from Noetherian rings to coherent rings, and in this process, finitely generated modules were replaced by finitely presented modules.  Recently, as  extending work of Stenstr\"{o}m's viewpoint, Gao and Wang \cite{GW} introduced the notion of weak injective modules. This class of modules was also investigated by Bravo,  Gillespie, and Hovey \cite{BGH} independently. In this process, finitely presented modules were  replaced by super finitely
presented modules (see \cite{GW0} or Section 2 for the definition). The fact  shows that  weak injective modules play a crucial role in the process of generalizing  homological properties from special rings to arbitrary
rings (see \cite{GH,GW,Zh1,Zh2,Zh3} for more details).

In 1965, Eilenberg and Moore first introduced the
theory of relative homological algebra in \cite{EM}. Since then the relative homological algebra,
especially the Gorenstein homological algebra, got a rapid
development. Nowadays, it has been developed to an advanced level (e.g. \cite{BGO1,Ch,CFH,EJ,Ho,Hu,HH,SSW}).
However, in  most results of Gorenstein homological algebra, the
condition `noetherian' is essential. In order to make
similar properties of Gorenstein homological algebra hold over general rings, one of methods is to define
a class of modules using an exact sequence of injective modules which is also exact after applying the covariant Hom functor with respect to weak injective modules.  But  these  modules are in fact stronger than the Gorenstein injective modules (see \cite[Def. 10.1.1]{EJ} for
the definition), and  weak injective modules are not contained in these modules in general. So our attempt is to find a new class of modules satisfying the following properties at the same time:
\begin{itemize}
  \item this new class of  modules are weaker than that of Gorenstein injective modules, and
  \item every weak injective module is Gorenstein weak injective.
\end{itemize}

 To solve this question,  we investigate a class of Gorenstein weak injective modules inspired by \cite{Gao}, that is, an $R$-module $M$ is called \emph{Gorenstein weak injective} if there exists an exact sequence of weak injective $R$-modules
$$
 \xymatrix@C=0.5cm{
  \cdots \ar[r] & E_1 \ar[r]^{} & E_0 \ar[r]^{} & E^0 \ar[r]^{} & E^1 \ar[r] & \cdots }
$$
such that $M=\operatorname{Coker}(E_1\rightarrow E_0)$ and the functor $\operatorname{Hom}_R(P,-)$ leaves this sequence exact whenever $P$ is a super finitely presented $R$-module with $pd_R(P)<\infty$.

This paper is organised as follows.
In Section 2, we first introduce the notions of Gorenstein weak injective and weak flat modules in terms of weak injective  and weak flat modules respectively,  and give some basic properties of them. Then we characterize rings over which all modules are Gorenstein weak injective and others. In Section 3, we  discuss a relation between weak cosyzygy and Gorenstein weak cosyzygy of a module, and investigate the stability of Gorenstein weak injective modules.

\section{Gorenstein weak injective and weak flat  modules}

In this section, we give the definitions of Gorenstein weak injective and weak flat modules and discuss some of properties of these modules. We first recall some terminologies and  preliminary. For more details, we refer the readers to \cite{EJ,GW0,GW}.

\begin{definition} {(\cite[Def. 6.1.1]{EJ})
Let $\mathcal {F}$ be a class of $R$-modules. By an
\emph{$\mathcal {F}$-preenvelope} of an $R$-module $M$, we mean a morphism
$\varphi: M \rightarrow F$ where $F\in \mathcal {F}$ such that for
any morphism $f: M\rightarrow F^{'}$ with $F^{'}\in \mathcal {F}$,
there exists a morphism $g:F\rightarrow F^{'}$ such that
$g\varphi=f$, that is, there is the following commutative diagram:
$$
\xymatrix{
  M \ar[dr]_{f} \ar[r]^{\varphi}
                & F \ar@{.>}[d]^{g}  \\
                & F'             }
$$
 If furthermore, when $F^{'}=F$ and $f=\varphi$,
the only such $g$ are automorphisms of $F$, then $\varphi: M
\rightarrow F$ is called an \emph{$\mathcal {F}$-envelope} of $M$.

Dually, one may give the notion of \emph{$\mathcal {F}$-(pre)cover} of an $R$-module.}
\end{definition}

We note that $\mathcal {F}$-envelopes and $\mathcal {F}$-covers may not
exist in general, but if they exist, they are unique up to
isomorphism.

We also note that  if the class $\mathcal {F}$ contains all injective $R$-modules, then
$\mathcal {F}$-preenvelopes are monic, and if the class $\mathcal {F}$ contains all projective $R$-modules, then
$\mathcal {F}$-precovers are epic.

In the process of generalizing  homological results  from special rings to arbitrary rings, the notion of super finitely presented modules plays a crucial role. Recall from \cite{GW0} that an $R$-module $M$ is called \emph{super finitely presented} if there exists an exact sequence $\cdots \rightarrow F_1\rightarrow F_0\rightarrow M\rightarrow 0$, where each $F_i$ is finitely generated and projective. Following this,  Gao and Wang \cite{GW} gave the notions of weak injective and weak flat modules in terms of super finitely presented modules, which are generalizations of  FP-injective and flat modules. Independently, these classes of modules were also considered by Bravo, Gillespie and Hovey \cite{BGH}, and in their paper, they used the notions of $\mbox{FP}_\infty$-injective and level modules rather than  weak injective  and weak flat modules.

\begin{definition}{(\cite[Def. 2.1]{GW})
An $R$-module $M$ is called \emph{weak injective} if $\mbox{Ext}^1_R(N,M)=0$ for any super finitely presented $R$-module $N$. A right $R$-module $M$ is called \emph{weak flat} if $\mbox{Tor}_1^R(M,N)=0$ for any super finitely presented $R$-module $N$.}
\end{definition}

The following implications are clear:
\begin{align*}
  \mbox{injective }R\mbox{-modules} & \Rightarrow\mbox{ FP-injective }R\mbox{-modules} \\
   & \Rightarrow\mbox{ weak injective }R\mbox{-modules}.\\
   \mbox{flat right }R\mbox{-modules} & \Rightarrow\mbox{ weak flat right }R\mbox{-modules}.
\end{align*}
Moreover, the class of FP-injective $R$-modules coincides with that of weak injective $R$-modules and the class of flat right $R$-modules coincides with that of weak flat right $R$-modules over a left coherent ring $R$ by \cite[Rem. 2.2]{GW}, and the class of injective $R$-modules coincides with that of FP-injective $R$-modules over a left Noetherian ring $R$ by \cite[Thm. 3]{Me}. Consequently, the class of injective $R$-modules coincides with that of weak injective $R$-modules over a Noetherian ring $R$.

We denote by $\mathcal{WI}(R)$ and $\mathcal{WF}(R^{\operatorname{op}})$ the classes of  weak injective and weak flat (right) $R$-modules respectively.  By  \cite[Thm. 3.4]{GH},  every $R$-module  has  a weak injective preenvelope. So for any $R$-module $M$, $M$ has a right $\mathcal{WI}(R)$-resolution.
Moreover, since every injective $R$-module is weak injective, every  right $\mathcal{WI}(R)$-resolution of $M$ is also exact, that is, there exists an exact sequence
$$\xymatrix@C=0.5cm{
  0 \ar[r]&  M\ar[r]& E^0\ar[r]& E^1\ar[r]&E^2\ar[r]& \cdots},$$
  where each $E^i$ is weak injective.

Moreover, following \cite{GH} or \cite[Thm. 12]{ET}, every right  $R$-module has a weak flat (pre)cover. So every right $R$-module $M$ has a left $\mathcal{WF}(R^{\operatorname{op}})$-resolution. Since every projective right $R$-module is weak flat, this resolution is also exact, that is, there is an exact sequence $$\cdots\rightarrow W_2\rightarrow W_1\rightarrow W_0\rightarrow M\rightarrow 0$$  with each $W_i$ weak flat.

\subsection{Gorenstein weak injective modules}

In homological algebra, finitely generated modules are well-behaved over Noetherian rings, and finitely presented modules are
 well-behaved over coherent rings. One of main reasons is that it always provides a resolution of a finitely generated  module by finitely generated projective modules.
 Thus in order to study homological properties in general rings, we first choose a class of modules that admit a resolution by finitely generated projective modules.
 We notice that a module $N$ is super finitely presented with
$pd_R(N)<\infty$ if and only if it has a finite resolution by finitely generated projective modules. Moreover, Gao and Wang \cite{Gao} used finitely presented modules of finite projective dimension to define Gorenstein $FP$-injective modules, and showed that they are well-behaved over coherent rings.  Thus using super finitely presented modules of finite projective dimension,
 we give the definition of Gorenstein weak injective modules as follows.

\begin{definition}\label{Gorenstein weak injective}
{ An $R$-module $M$ is called \emph{Gorenstein weak injective} if
there exists an exact sequence of weak injective $R$-modules
$$
 \mathbb{W}=\xymatrix@C=0.5cm{
  \cdots \ar[r] & W_1 \ar[r]^{} & W_0 \ar[r]^{} & W^0 \ar[r]^{} & W^1 \ar[r] & \cdots }
$$
such that $M=\operatorname{Coker}(W_1\rightarrow W_0)$ and the
functor $\operatorname{Hom}_R(N,-)$ leaves this sequence exact
whenever $N$ is a super finitely presented $R$-module with
$pd_R(N)<\infty$.}
\end{definition}

We denote by $\mathcal{GWI}(R)$ the class of
Gorenstein weak injective $R$-modules.

\begin{remark}\label{remark1}
\begin{itemize}
\item[]
  \item [(1)] By definition,  every weak injective $R$-module is Gorenstein weak injective.
  \item[(2)] Since every FP-injective $R$-module is weak injective, every
Gorenstein FP-injective $R$-module (in the sense of Gao's
definition, see \cite{Gao}) is Gorenstein weak injective. If $R$ is a left coherent
ring, then the class of  Gorenstein weak injective $R$-modules
coincides with the class of  Gorenstein FP-injective $R$-modules.
Moreover, we have the following implications by \cite[Prop.
2.5]{Gao}:
\begin{align*}
    \mbox{Gorenstein injective }R\mbox{-modules}&\Rightarrow\mbox{ Gorenstein FP-injective }R\mbox{-modules} \\
    &\Rightarrow\mbox{ Gorenstein weak injective }R\mbox{-modules}.
\end{align*}
If $R$ is an $n$-Gorenstein ring (i.e. a left and right Noetherian
ring with self-injective dimension at most $n$ on both sides for
some non-negative integer $n$), then these three kinds of
$R$-modules coincide.
  \item [(3)] The class of Gorenstein weak injective $R$-modules is closed
under direct products  by definition.
  \item [(4)] If $
 \mathbb{W}=\xymatrix@C=0.5cm{
  \cdots \ar[r] & W_1 \ar[r]^{} & W_0 \ar[r]^{} & W^0 \ar[r]^{} & W^1 \ar[r] & \cdots }
$ is an exact sequence of weak injective $R$-modules such that the
functor $\operatorname{Hom}_R(N,-)$ leaves this sequence exact
whenever $N$ is a super finitely presented $R$-module with
$pd_R(N)<\infty$, then by symmetry, all the images, the kernels and
the cokernels of $\mathbb{W}$ are Gorenstein weak injective.
\end{itemize}
\end{remark}

\begin{proposition}\label{functorchar}
Let $M$ be a Gorenstein weak injective $R$-module. Then we have
$\operatorname{Ext}^i_R(N,M)=0$  whenever $N$ is a super finitely
presented $R$-module with $pd_R(N)<\infty$ and $i\geq 1$.
\end{proposition}

\begin{proof}
Since $M$ is a Gorenstein weak injective $R$-module, by Definition $\ref{Gorenstein weak injective}$,
there exists an exact sequence $0\rightarrow M\rightarrow
W^0\rightarrow M^1\rightarrow 0$ with $W^0$ weak injective and $M^1$
Gorenstein weak injective, such that the functor
$\operatorname{Hom}_R(N,-)$ leaves this sequence exact whenever $N$
is a super finitely presented $R$-module with $pd_R(N)<\infty$.
Moreover, consider the following exact sequence
$$
\xymatrix@C=0.4cm{
  0 \ar[r] & \mbox{Hom}_R(N,M) \ar[r]^{} & \mbox{Hom}_R(N,W^0) \ar[r]^{} & \mbox{Hom}_R(N,M^1) \ar[r]^{} & \mbox{Ext}^1_R(N,M) \ar[r]^{} &  \mbox{Ext}^1_R(N,W^0) }.
$$
Since $W^0$ is  weak injective, we have $\mbox{Ext}^1_R(N,W^0)=0$,
and hence $\mbox{Ext}^1_R(N,M)=0$. Moreover, since $M^1$ is
Gorenstein weak injective, we also have $\mbox{Ext}^1_R(N,M^1)=0$.
Consider the following exact sequence
$$
\xymatrix@C=0.5cm{
  0=\mbox{Ext}^1_R(N,W^0) \ar[r] & \mbox{Ext}^1_R(N,M^1) \ar[r]^{} & \mbox{Ext}^2_R(N,M) \ar[r]^{} & \mbox{Ext}^2_R(N,W^0)}.
$$
Note that $\mbox{Ext}^2_R(N,W^0)=0$ by \cite[Prop. 3.1]{GW}, and
hence   $\mbox{Ext}^2_R(N,M)\cong \mbox{Ext}^1_R(N,M^1)=0$.
  We repeat the argument by replacing $M^1$ with $M$ to get a Gorenstein weak injective $R$-module $M^2$ and the isomorphisms $\mbox{Ext}^3_R(N,M)\cong \mbox{Ext}^2_R(N,M^1)\cong\mbox{Ext}^1_R(N,M^2)=0$. Continuing this process, we may obtain a Gorenstein weak injective $R$-module $M^{i-1}$ and the isomorphisms $\mbox{Ext}^i_R(N,M)\cong \mbox{Ext}^{i-1}_R(N,M^1)\cong\cdots\cong\mbox{Ext}^1_R(N,M^{i-1})=0$ for any $i\geq 1$.
\end{proof}

The following proposition shows that we may simplify the definition
of Gorenstein weak injective $R$-modules.

\begin{proposition}\label{simplify}
The following are equivalent for an $R$-module $M$:
\begin{itemize}
  \item [(1)] $M$ is Gorenstein weak injective;
  \item [(2)]  There
exists an exact sequence of weak injective $R$-modules
$$
 \mathbb{W}=\xymatrix@C=0.5cm{
  \cdots \ar[r] & W_1 \ar[r]^{} & W_0 \ar[r]^{} & W^0 \ar[r]^{} & W^1 \ar[r] & \cdots }
$$
such that $M=\operatorname{Coker}(W_1\rightarrow W_0)$;
  \item[(3)] There exists an exact sequence $\cdots \rightarrow W_1\rightarrow
W_0\rightarrow M\rightarrow 0$, where each $W_i$ is weak injective;
  \item  [(4)] There exists an exact sequence $0 \rightarrow L\rightarrow
W\rightarrow M\rightarrow 0$, where $W$ is weak injective and $L$ is
Gorenstein weak injective.
\end{itemize}
\end{proposition}

\begin{proof}
(1) $\Rightarrow$ (4) $\Rightarrow$ (3) are trivial.

(3) $\Rightarrow$ (2). Since every $R$-module has a weak injective
preenvelope by \cite[Thm. 3.4]{GH}, we may easily get an exact
sequence $0\rightarrow M \rightarrow W^0\rightarrow
W^1\rightarrow\cdots$, where each $W^i$ is weak injective.
Assembling this sequence with the sequence given in (3), we get the
following exact sequence
$$
 \mathbb{W}=\xymatrix@C=0.5cm{
  \cdots \ar[r] & W_1 \ar[r]^{} & W_0 \ar[r]^{} & W^0 \ar[r]^{} & W^1 \ar[r] & \cdots }
$$
such that $M=\operatorname{Coker}(W_1\rightarrow W_0)$.

(2) $\Rightarrow$(1).  By Definition $\ref{Gorenstein weak injective}$,
it suffices to show that the complex
$\operatorname{Hom}_R(N,\mathbb{W})$ is exact whenever $N$ is
a super finitely presented $R$-module with $pd_R(N)<\infty$.

We use induction on $n=pd_R(N)<\infty$. The case $n=0$ is trivial.
Let $n\geq 1$, and assume that the result holds for the case $n-1$.
Consider an exact sequence $0\rightarrow K\rightarrow P_0\rightarrow
N\rightarrow 0$, where $P_0$ is finitely generated projective and
$K$ is super finitely presented. Then $pd_R(K)=n-1$. Since each term
of $\mathbb{W}$ is weak injective, we may get the following
exact sequence of complexes
$$
0\rightarrow
\mbox{Hom}_R(N,\mathbb{W})\rightarrow\mbox{Hom}_R(P_0,{\mathbb{W}})\rightarrow
\mbox{Hom}_R(K,\mathbb{W})\rightarrow 0.
$$
Clearly, the complex $\mbox{Hom}_R(P_0,\mathbb{W})$ is exact.
Moreover, the complex $\mbox{Hom}_R(K,\mathbb{W})$ is also
exact by the induction hypothesis. So the complex
$\mbox{Hom}_R(N,\mathbb{W})$ is exact, and hence $M$ is
Gorenstein weak injective.
\end{proof}

Following this and \cite[Prop. 2.6]{GW}, we have

\begin{corollary}
Let $I$ be a directed set, and $\{M_i\}_{i\in I}$ a direct system of $R$-modules. If every $M_i$ is
Gorenstein weak injective, then  the direct limit $\underrightarrow{\lim}  M_i$ is
Gorenstein weak injective.
\end{corollary}

\begin{proposition}
Given an exact sequence $0\rightarrow L\rightarrow M\rightarrow
N\rightarrow 0$ of $R$-modules.
\begin{itemize}
  \item [(1)] If $M$ is weak injective and $L$ is Gorenstein weak
injective, then $N$ is Gorenstein weak
injective.
  \item [(2)] If $L$ is weak injective and  $N$ is Gorenstein weak
injective, then $M$  is Gorenstein weak
injective.
\end{itemize}
\end{proposition}

\begin{proof}
(1) follows from Proposition \ref{simplify}.

 (2) Assume that $N$ is Gorenstein weak injective. Then, by Proposition \ref{simplify},
 there exists an exact sequence $0\rightarrow K\rightarrow W\rightarrow N\rightarrow 0$,
 where $K$ is Gorenstein weak injective and $W$ is weak injective. Consider the following pull-back diagram:
\begin{gather}\label{1}
\begin{split}
\xymatrix@=20pt{
&&0\ar[d]&0\ar[d]&\\
&&K\ar[d]\ar@{=}[r]&K\ar[d]&\\
0\ar[r]&L\ar[r]\ar@{=}[d]&W'\ar[r]\ar[d]&W\ar[r]\ar[d]&0\\
0\ar[r]&L\ar[r]&M\ar[r]\ar[d]&N\ar[r]\ar[d]&0\\
&&0&0& }
\end{split}
\end{gather}
Since $L$ and $W$ are weak injective, it follows from the middle row
in the diagram ($\ref{1}$) that $W'$ is weak injective. Moreover, by
the middle column in the  diagram ($\ref{1}$), we have that $M$ is Gorenstein weak injective.
\end{proof}

\begin{definition}\label{Gorenstein weak dimension}
{ The \emph{Gorenstein weak injective dimension} of an $R$-module $M$,
denoted by $Gwid_R(M)$, is defined to be $\operatorname{inf}\{n\mid
\mbox{there is an exact sequence }0\rightarrow M\rightarrow
G^0\rightarrow G^1\rightarrow \cdots\rightarrow G^n\rightarrow 0
\mbox{ with }  G^i \mbox{ Gorenstein weak injective for any } 0\leq
i\leq n\}.$ If no such $n$ exists, set $Gwid_R(M)=\infty$.}
\end{definition}

\subsection{Gorenstein weak flat modules}

Now we give the definition of Gorenstein weak flat modules.

\begin{definition}\label{Gorenstein weak flat}
{A right $R$-module $M$ is called \emph{Gorenstein weak flat} if
there exists an exact sequence of weak flat right $R$-modules
$$
 \mathbb{W}=\xymatrix@C=0.5cm{
  \cdots \ar[r] & W_1 \ar[r]^{} & W_0 \ar[r]^{} & W^0 \ar[r]^{} & W^1 \ar[r] & \cdots }
$$
such that $M=\operatorname{Coker}(W_1\rightarrow W_0)$ and the
functor $-\otimes_RN$ leaves this sequence exact
whenever $N$ is a super finitely presented $R$-module with
$pd_R(N)<\infty$.}
\end{definition}

We denote by $\mathcal{GWF}(R)$ the class of
Gorenstein weak flat $R$-modules.

\begin{remark}\label{remark2}
\begin{itemize}
\item[]
  \item [(1)] By definition,  every weak flat right $R$-module is Gorenstein weak flat.
  \item [(2)] The class of Gorenstein weak flat right $R$-modules is closed
under  direct  sums and direct products by definition and \cite[Thm. 2.13]{GW}.
  \item [(3)] If $
 \mathbb{W}=\xymatrix@C=0.5cm{
  \cdots \ar[r] & W_1 \ar[r]^{} & W_0 \ar[r]^{} & W^0 \ar[r]^{} & W^1 \ar[r] & \cdots }
$ is an exact sequence of weak flat right $R$-modules such that the
functor $-\otimes_RN$ leaves this sequence exact
whenever $N$ is a super finitely presented $R$-module with
$pd_R(N)<\infty$, then by symmetry, all the images, the kernels and
the cokernels of $\mathbb{W}$ are Gorenstein weak flat.
\end{itemize}
\end{remark}

As a similar argument to that of Proposition \ref{simplify}, we have

\begin{proposition}\label{simplify2}
The following are equivalent for a right $R$-module $M$:
\begin{itemize}
  \item [(1)] $M$ is Gorenstein weak flat;
  \item [(2)]  There
exists an exact sequence of weak flat right $R$-modules
$$
 \mathbb{W}=\xymatrix@C=0.5cm{
  \cdots \ar[r] & W_1 \ar[r]^{} & W_0 \ar[r]^{} & W^0 \ar[r]^{} & W^1 \ar[r] & \cdots }
$$
such that $M=\operatorname{Coker}(W_1\rightarrow W_0)$;
  \item [(3)] There exists an exact sequence $0 \rightarrow M \rightarrow W^0\rightarrow
W^1\rightarrow \cdots$, where each $W^i$ is weak flat;
  \item [(4)] There exists an exact sequence $0 \rightarrow M\rightarrow
W\rightarrow N\rightarrow 0$, where $W$ is weak flat and $N$ is
Gorenstein weak flat.
\end{itemize}
\end{proposition}

Recall from \cite[Def. 10.3.1]{EJ} that  a right $R$-module $M$ is called \emph{ Gorenstein  flat} if
there exists an exact sequence of  flat right $R$-modules
$$
 \mathbb{F}=\xymatrix@C=0.5cm{
  \cdots \ar[r] & F_1 \ar[r]^{} & F_0 \ar[r]^{} & F^0 \ar[r]^{} & F^1 \ar[r] & \cdots }
$$
such that $M=\operatorname{Coker}(F_1\rightarrow F_0)$ and the
functor $-\otimes_RI$ leaves this sequence exact
whenever $I$ is an injective $R$-module.

Recall that a ring $R$ is called \emph{$n$-FC} if it is left and right coherent and $FP\mbox{-}id_R({_RR})\leq n$ and $FP\mbox{-}id_R({R_R})\leq n$, where the symbol $FP\mbox{-}id_R(-)$ denotes the FP-injective dimension of modules. The following proposition shows that the class of Gorenstein weak flat modules is larger than that of Gorenstein  flat modules, and they have no difference over $n$-FC rings.
\begin{proposition}
Every Gorenstein flat right $R$-module is  Gorenstein weak flat. The converse holds if $R$ is an $n$-FC ring.
\end{proposition}

\begin{proof}
Let $M$ be a Gorenstein flat right $R$-module. Then, by definition, there is an exact sequence of flat right $R$-modules
$$
 \mathbb{F}=\xymatrix@C=0.5cm{
  \cdots \ar[r] & F_1 \ar[r]^{} & F_0 \ar[r]^{} & F^0 \ar[r]^{} & F^1 \ar[r] & \cdots }
$$
with $M=\mbox{Coker}(F_1\rightarrow F_0)$. Note that each flat right $R$-module is weak flat. It follows from Proposition \ref{simplify2} that $M$ is Gorenstein weak flat.

Conversely, assume that $M$ is a Gorenstein weak flat right $R$-module. Then by Definition \ref{Gorenstein weak flat}, there
is an exact sequence of weak flat right $R$-modules
$$
 \mathbb{W}=\xymatrix@C=0.5cm{
  \cdots \ar[r] & W_1 \ar[r]^{} & W_0 \ar[r]^{} & W^0 \ar[r]^{} & W^1 \ar[r] & \cdots }
$$
such that $M=\operatorname{Coker}(W_1\rightarrow W_0)$. Note that, over a left coherent ring, the class of weak flat right modules coincides with that of flat right modules. So $\mathbb{W}$ is also an exact sequence of flat right $R$-modules. We will argue that the complex $\mathbb{W}\otimes_R E$ is exact for any left $R$-module $E$ with $FP\mbox{-id}_R(E)<\infty$. Since $R$ is an $n$-FC ring, then we have that $FP\mbox{-}id_R(E)<\infty$ if and only if  $fd_R(E)\leq n$ by \cite[Prop. 3.6]{DC93}. So it suffices to show that $\mathbb{W}\otimes_R E$ is exact for any left $R$-module $E$ with $fd_R(E)\leq n$. We use induction on $n \geq fd_R(E)$. The case $n=0$ is trivial. Consider an exact sequence $0\rightarrow L\rightarrow P\rightarrow E\rightarrow 0$ with $P$ projective. Then we have the following exact sequence of complex
$$
0\longrightarrow \mathbb{W}\otimes_RL\longrightarrow \mathbb{W}\otimes_RP\longrightarrow \mathbb{W}\otimes_RE\longrightarrow 0
$$
Note that since $fd_R(E)\leq n$, we have $fd_R(L)\leq n-1$. By the hypothesis induction, $\mathbb{W}\otimes_RL$ is exact. Moreover, $\mathbb{W}\otimes_RP$ is exact, so $\mathbb{W}\otimes_RE$ is also exact. In particular, $\mathbb{W}\otimes_R I$ is exact for any injective  $R$-module $I$. Therefore, $M$ is Gorenstein flat.
\end{proof}

Following \cite[Thm. 2.10 and Props. 2.11, 2.12]{GW}, we can easily get

\begin{proposition}\label{char}
\begin{itemize}
\item[]
  \item [(1)] If a right $R$-module $M$ is
Gorenstein weak flat, then $M^+$ is
Gorenstein weak injective.
  \item [(2)] If an $R$-module $M$ is
Gorenstein weak injective, then $M^+$ is
Gorenstein weak flat.
  \item [(3)] If an $R$-module $M$ is
Gorenstein weak injective, then $M^{++}$ is
Gorenstein weak injective.
  \item [(4)] If a right $R$-module $M$ is
Gorenstein weak flat, then $M^{++}$ is
Gorenstein weak flat.
\end{itemize}
\end{proposition}

\begin{definition}
{The \emph{Gorenstein weak flat dimension} of an $R$-module $M$,
denoted by $Gwfd_R(M)$, is defined to be $\operatorname{inf}\{n\mid
\mbox{there is an exact sequence }0\rightarrow G^n\rightarrow
\cdots\rightarrow G^1\rightarrow G^0\rightarrow M\rightarrow 0
\mbox{ with }  G^i \mbox{ Gorenstein weak flat for any } 0\leq
i\leq n\}.$ If no such $n$ exists, set $Gwfd_R(M)=\infty$.}
\end{definition}

\subsection{Rings over which every module is
Gorenstein weak injective}

We now give a characterization for  rings whose every module is
Gorenstein weak injective and, meanwhile, every right module is
Gorenstein weak flat as follows.

\begin{proposition}\label{characterization}
The following are equivalent:
\begin{itemize}
  \item [(1)] Every $R$-module is Gorenstein weak injective;
  \item [(2)] Every right $R$-module is Gorenstein weak flat;
  \item [(3)] Every projective $R$-module is weak injective;
  \item [(4)] Every flat $R$-module is weak injective;
  \item [(5)] Every injective right $R$-module is weak flat;
  \item [(6)] $R$ is weak injective as an $R$-module.
\end{itemize}
\end{proposition}

\begin{proof}
(1) $\Rightarrow$ (3). Let $P$ be a projective $R$-module. Then $P$
is Gorenstein weak injective by hypothesis. So there exists an exact
sequence $0\rightarrow K\rightarrow W\rightarrow P\rightarrow 0$,
where $W$ is weak injective. Since $P$ is projective, this sequence is
split, and hence $P$ is weak injective as a direct summand of $W$ by
\cite[Prop. 2.3]{GW}.

(3) $\Rightarrow$ (1). Let $M$ be any $R$-module. If every
projective $R$-module is weak injective, then by assembling a
projective resolution of $M$ with its weak injective resolution, we
may get the following exact sequence of weak injective $R$-modules
$$
\xymatrix@C=0.5cm{
  \cdots \ar[r] & P_1 \ar[r]^{} & P_0 \ar[r]^{} & W^0 \ar[r]^{} & W^1 \ar[r] & \cdots }
$$
such that $M=\operatorname{Ker}(W^0\rightarrow W^1)$. Thus $M$ is
Gorenstein weak injective by Proposition $\ref{simplify}$.

(2) $\Rightarrow$ (5). Let $I$ be an injective right $R$-module. Then $I$
is Gorenstein weak flat by hypothesis. So there exists an exact
sequence $0\rightarrow I\rightarrow W\rightarrow N\rightarrow 0$,
where $W$ is weak flat. Since $I$ is injective, this sequence is
split, and hence $I$ is weak flat as a direct summand of $W$ by
\cite[Prop. 2.3]{GW}.

(5) $\Rightarrow$ (2). Let $M$ be any $R$-module. If every
injective right $R$-module is weak flat, then by assembling an
injective resolution of $M$ with its weak flat resolution, we
may get the following exact sequence of weak flat right $R$-modules
$$
\xymatrix@C=0.5cm{
  \cdots \ar[r] & W_1 \ar[r]^{} & W_0 \ar[r]^{} & I^0 \ar[r]^{} & I^1 \ar[r] & \cdots }
$$
such that $M=\operatorname{Coker}(W_1\rightarrow W_0)$. Thus $M$ is
Gorenstein weak flat by Proposition $\ref{simplify2}$.

(4) $\Leftrightarrow$ (5) $\Leftrightarrow$ (6) The proof is similar to that of \cite[Prop. 2.17]{GW}. (4) $\Rightarrow$ (3) $\Rightarrow$ (6) are trivial.
\end{proof}

\begin{corollary}
Let $R$ be a left Noetherian ring. Then the following are equivalent:

$(1)$ $R$ is quasi-Frobenius;

$(2)$ Every $R$-module is Gorenstein weak injective;

$(3)$ Every right $R$-module is Gorenstein weak flat.
\end{corollary}

\begin{proof}
(1) $\Rightarrow$ (2). It is obvious, since every $R$-module is Gorenstein injective over quasi-Frobenius rings by \cite[Prop. 2.6]{BM}.

(2) $\Leftrightarrow$ (3) follow from Proposition $\ref{characterization}$.

(3) $\Rightarrow$ (1) follows from Proposition
$\ref{characterization}$ and the fact that the injective $R$-modules
coincide with the weak injective $R$-modules over a left Noetherian
ring $R$.
\end{proof}

In \cite{GW}, Gao and Wang defined the left super finitely presented dimension of a ring R to be
$l.sp.gldim(R)=\sup\{pd_RM\mid M \mbox{ is a super finitely presented  $R$-module}\}$.
By \cite[Thm. 3.8]{GW},
\begin{align*}
l.sp.gldim(R) &=\mbox{sup}\{wid_R(M)\mid M \mbox{ is any $R$-module}\}  \\
  & =\mbox{sup}\{wfd_R(M)\mid M \mbox{ is any right $R$-module}\}.
                                                                     \end{align*}

\begin{proposition} If $l.sp.gldim(R)<\infty$, then we have
\begin{itemize}
  \item [(1)] Every Gorenstein weak injective $R$-module is weak injective;
  \item [(2)] Every Gorenstein weak flat right $R$-module is weak flat.
\end{itemize}
\end{proposition}

\begin{proof}
(1) Assume that $l.sp.gldim(R)=n<\infty$ and let $M$ be a Gorenstein weak
injective $R$-module. The case $n=0$ is trivial. Let $n\geq 1$.
Since $M$ is Gorenstein weak injective, there is an exact sequence
$$
\xymatrix@C=0.5cm{
  \cdots \ar[r] & W_1 \ar[r]^{} & W_0 \ar[r]^{} & M \ar[r]^{} & 0}
$$
with each $W_i$  weak injective. Let
$K_n=\mbox{Ker}(W_{n-1}\rightarrow W_{n-2})$. Then we get an exact
sequence
$$
\xymatrix@C=0.5cm{
  0\ar[r]&K_n\ar[r]& W_{n-1}\ar[r]&\cdots \ar[r] & W_1 \ar[r]^{} & W_0 \ar[r]^{} & M \ar[r]^{} & 0}.
$$
By hypothesis, $wid_R(K_n)\leq n$, and hence $M$ is weak injective
by \cite[Prop. 3.3]{GW}.

(2) The proof is similar to that of (1).
\end{proof}

The above proposition shows that the class of weak injective
$R$-modules coincides with that of Gorenstein  weak injective
$R$-modules, and  the class of weak flat right
$R$-modules coincides with that of Gorenstein  weak flat right
$R$-modules  over a ring $R$ satisfying $l.sp.gldim(R)<\infty$.

The next proposition also gives a description of  rings over which
all Gorenstein weak injective  $R$-modules are weak injective from
the viewpoint of Gorenstein weak injective dimension of modules.

\begin{proposition}
The following are equivalent:
\begin{itemize}
  \item [(1)] Every Gorenstein weak injective  $R$-module is weak injective;
  \item [(2)] For any $R$-module $M$, $Gwid_R(M)=wid_R(M)$.
\end{itemize}
\end{proposition}

\begin{proof}
(1) $\Rightarrow$ (2). Let $M$ be an $R$-module. Since every weak
injective $R$-module is Gorenstein weak injective, it is obvious
that $Gwid_R(M)\leq wid_R(M)$. Thus it suffices to show that
$wid_R(M)\leq Gwid_R(M)$. Without loss of generality, we assume that
$Gwid_R(M)=n<\infty$ for some non-negative integer $n$. Then, by
Definition $\ref{Gorenstein weak dimension}$, there is an exact
sequence $0\rightarrow M\rightarrow G^0\rightarrow G^1\rightarrow
\cdots\rightarrow G^n\rightarrow 0$ with   $G^i$ Gorenstein weak
injective for any  $0\leq i\leq n$. Note that each $G^i$ is weak
injective by hypothesis. Thus, $wid_R(M)\leq n=Gwid_R(M)$ by
\cite[Prop. 3.3]{GW}, as desired.

(2) $\Rightarrow$ (1)  is trivial.
\end{proof}

Similarly, we have

\begin{proposition}
The following are equivalent:
\begin{itemize}
  \item [(1)] Every Gorenstein weak flat   $R$-module is weak flat;
  \item [(2)] For any $R$-module $M$, $Gwfd_R(M)=wfd_R(M)$.
\end{itemize}
\end{proposition}

\section{The cosyzygy and stability of Gorenstein weak injective modules}

Let $\mathcal{X}$ and $\mathcal{Y}$ be two classes of $R$-modules. We write
$\mathcal{X}\bot\mathcal{Y}$ if $\mbox{Ext}_R^1(X,Y)=0$ for any $X\in\mathcal{X}$ and $Y\in\mathcal{Y}$.

\begin{remark}
\begin{itemize}
  \item [(1)] If $\mathcal{WI}(R)\bot\mathcal{GWI}(R)$, then $\mathcal{GWI}(R)$ is closed under extensions. Indeed, let
$0\rightarrow L\rightarrow M\rightarrow N\rightarrow 0$ be an exact sequence of $R$-modules
with $L,N\in\mathcal{GWI}(R)$.
By Proposition \ref{simplify}, there exist exact  sequences
\begin{gather*}
    \xymatrix@C=15pt{
\cdots\ar[r] &E^{'}_1 \ar[r]^{d^{'}_1} & E^{'}_0 \ar[r]^{d^{'}_0} &L \ar[r]
& 0}\ {\rm and} \ \
\xymatrix@C=15pt{
\cdots\ar[r] &E^{''}_1 \ar[r]^{d^{''}_1} & E^{''}_0 \ar[r]^{d^{''}_0} &N \ar[r]
& 0}
\end{gather*}
with all $E^{'}_i$, $E^{''}_i$ in $\mathcal{WI}(R)$ and all $\mbox{Ker}d^{'}_i$, $\mbox{Ker}d^{''}_i$ in $\mathcal{GWI}(R)$.
 Consider the following diagram
$$\xymatrix@R=15pt@C=20pt{
&E^{'}_0\ar[d]^{d^{'}_0}&&E^{''}_0\ar[d]^{d^{''}_0}&\\
0\ar[r]&L\ar[r]^f&M\ar[r]^g&N\ar[r]&0.}
$$
Since $\mbox{Ext}^1_R(E^{''}_0,L)=0$, we get an epimorphism
$\mbox{Hom}_R(E^{''}_0,M) \twoheadrightarrow \mbox{Hom}_R(E^{''}_0,N)$ and
there exists $\alpha:E^{''}_0\rightarrow M$ such that $d^{''}_0=g\alpha$.
Putting  $E_0:=E^{'}_0\oplus E^{''}_0$ and $d_0:=(fd^{'}_0 \ \alpha)$,
then we obtain the following commutative diagram with exact columns and rows
$$\xymatrix@R=15pt@C=20pt{
0\ar[r]&E^{'}_0\ar[r]^{\!\!{1\choose 0}}\ar[d]^{d^{'}_0}&E_0\ar[r]^{\ (0\ 1)}\ar[d]^{d_0}&E^{''}_0\ar[r]\ar[d]^{d^{''}_0}&0\\
0\ar[r]&L\ar[r]^f\ar[d]&M\ar[r]^g\ar[d]&N\ar[r]\ar[d]&0\\
&0&0&0.&}
$$
Repeating this process, we may get an exact  sequence
$$\xymatrix@C=0.5cm{
\cdots \ar[r] &E_1\ar[r] & E_0 \ar[r]^{} &M \ar[r]& 0}$$  with all $E_i=E^{'}_i\oplus E^{''}_i$ weak injective. Thus
$M\in\mathcal{GWI}(R)$ by Proposition  \ref{simplify}.
  \item [(2)] If $R$ is a Gorenstein ring, then $\mathcal{GWI}(R)$ is closed under extensions.
  Indeed, in this case, $\mathcal{GWI}(R)$ coincides with the class of Gorenstein injective $R$-modules.
\end{itemize}
\end{remark}

In the following,  we always assume that the ground ring is a ring $R$ over which the class of Gorenstein
weak injective $R$-modules is closed under extensions.

Consider the following exact sequence
$$
\xymatrix@C=0.5cm{
  0 \ar[r] & M \ar[r]^{d^0} & W^0 \ar[r]^{d^1} & W^1 \ar[r]^{} & \cdots},
$$
where each $W^i$ is weak injective. Let $V^i=\mbox{Coker}d^{i-1}$
for any $i\geq 1$. Then we call $V^i$ an $i$th \emph{weak
cosyzygy} of $M$.

Similarly, if each $W^i$ is Gorenstein weak injective in the above
sequence, then  we call $V^i$ an $i$th \emph{Gorenstein weak
cosyzygy} of $M$.

We will investigate a relation between weak cosyzygy and
Gorenstein weak cosyzygy of an $R$-module as follows.

Since every weak injective $R$-module is Gorenstein weak injective,
it is obvious that every $i$th weak cosyzygy of  an $R$-module $M$
is an $i$th Gorenstein weak cosyzygy of $M$. The following theorem
shows that the converse holds in some cases.

\begin{theorem}\label{thmcosyzygy}
Let  $n$ be a positive integer and $V^n$
an $n$th Gorenstein weak cosyzygy of an $R$-module $M$. Then $V^n$
is an $n$th weak cosyzygy of some $R$-module $N$, and there is an
exact sequence $0\rightarrow G\rightarrow N\rightarrow M\rightarrow
0$, where $G$ is Gorenstein weak injective.
\end{theorem}

\begin{proof}
We use induction on $n$. For the case $n=1$, there is an exact
sequence $0\rightarrow M\rightarrow G^0\rightarrow V^1\rightarrow 0$
with $G^0$ Gorenstein weak injective. Moreover, there is an exact
sequence $0\rightarrow G\rightarrow W^0 \rightarrow G^0\rightarrow
0$  with $W^0$ weak injective and $G$ Gorenstein weak injective.

Consider the following pull-back diagram:
\begin{gather}\label{2}
\begin{split}
\xymatrix@=20pt{
&0\ar[d]&0\ar[d]&&\\
&G\ar@{=}[r]\ar[d]&G\ar[d]&&\\
0\ar[r]&N\ar[r]\ar[d]&W^0\ar[r]\ar[d]&V^1\ar@{=}[d]\ar[r]&0\\
0\ar[r]&M\ar[r]\ar[d]&G^0\ar[r]\ar[d]&V^1\ar[r]&0\\
&0&0&& }
\end{split}
\end{gather}
It follows from the middle row in the  diagram ($\ref{2}$) that
$V^1$ is a $1$st weak cosyzygy of an $R$-module $N$. Moreover, we
get the desired exact sequence $0\rightarrow G\rightarrow
N\rightarrow M\rightarrow 0$ from the second column in the  diagram
($\ref{2}$).

Now let $n\geq 2$ and suppose that the result holds for the case
$n-1$. Let $V^n$ be an $n$th Gorenstein weak cosyzygy of $M$. Then
we have the following  exact sequence
$$
\xymatrix@C=0.5cm{
  0\ar[r]&M\ar[r]& G^{0}\ar[r]&G^1\ar[r]&\cdots \ar[r] & G^{n-1} \ar[r]&V^n \ar[r]^{} & 0},
$$
where each $G^i$ is Gorenstein weak injective. Since $G^{n-1}$ is
Gorenstein weak injective, there is an exact $0\rightarrow
G'\rightarrow W^{n-1}\rightarrow G^{n-1}\rightarrow 0$ with $G'$
Gorenstein weak injective and $W^{n-1}$ weak injective.

Consider the following pull-back diagrams:
\begin{gather}\label{3}
\begin{split}
\xymatrix@=20pt{
&0\ar[d]&0\ar[d]&&\\
&G'\ar@{=}[r]\ar[d]&G'\ar[d]&&\\
0\ar[r]&V'\ar[r]\ar[d]&W^{n-1}\ar[r]\ar[d]&V^n\ar@{=}[d]\ar[r]&0\\
0\ar[r]&V^{n-1}\ar[r]\ar[d]&G^{n-1}\ar[r]\ar[d]&V^n\ar[r]&0\\
&0&0&& }
\end{split}
\end{gather}
and
\begin{gather}\label{4}
\begin{split}
\xymatrix@=20pt{
&&0\ar[d]&0\ar[d]&\\
&&G'\ar@{=}[r]\ar[d]&G'\ar[d]&\\
0\ar[r]&V^{n-2}\ar[r]\ar@{=}[d]&G''\ar[r]\ar[d]&V'\ar[d]\ar[r]&0\\
0\ar[r]&V^{n-2}\ar[r]&G^{n-2}\ar[r]\ar[d]&V^{n-1}\ar[r]\ar[d]&0\\
&&0&0& }
\end{split}
\end{gather}
where $V^{n-i}=\mbox{Coker}(G^{n-i-2}\rightarrow G^{n-i-1})$ for
$i=1,2$. For the exact sequence $0\rightarrow G'\rightarrow
G''\rightarrow G^{n-2}\rightarrow 0$ in the diagram ($\ref{4}$),
since $G'$ and $G^{n-2}$ are Gorenstein weak injective, $G''$ is
also Gorenstein weak injective. Moreover, by
 the middle row in the diagram ($\ref{4}$), we have that $V'$ is an $(n-1)$st Gorenstein weak
 cosyzygy of $M$. Thus $V'$ is an $(n-1)$st weak cosyzygy of some $R$-module $N$ by the
 induction hypothesis. In addition, by assembling the middle row in the diagram ($\ref{3}$),
 we may get that $V^n$ is an $n$th weak cosyzygy of $N$, as desired.
\end{proof}

In the following, we will consider the stability of Gorenstein weak
injective $R$-modules, which  shows that an iteration of the procedure used to describe
the class of Gorenstein  weak injective modules yields exactly the class of
Gorenstein  weak injective modules.

We begin with the following question, which is inspired by \cite{SSW} but different.

\begin{question} {\rm Given an exact sequence of Gorenstein weak injective $R$-modules
$$
\mathbb{G}=\xymatrix@C=0.5cm{
  \cdots \ar[r] & G_1 \ar[r]^{} & G_0 \ar[r]^{} & G^0 \ar[r]^{} & G^1 \ar[r] & \cdots }
$$
such that $M=\mbox{Coker}(G_1\rightarrow G_0)$ and the functor
$\operatorname{Hom}_R(N,-)$ leaves this sequence exact whenever $N$
is a super finitely presented $R$-module with $pd_R(N)<\infty$, is
$M$ Gorenstein weak injective?}
\end{question}

As a similar argument to the proof of $(2)\Rightarrow (1)$ in Proposition \ref{simplify}
and using Proposition \ref{functorchar}, the above question is equivalent to the following

\begin{question} {\rm Given an exact sequence of Gorenstein weak injective $R$-modules
$$
\mathbb{G}=\xymatrix@C=0.5cm{
  \cdots \ar[r] & G_1 \ar[r]^{} & G_0 \ar[r]^{} & G^0 \ar[r]^{} & G^1 \ar[r] & \cdots }
$$
such that $M=\mbox{Coker}(G_1\rightarrow G_0)$, is
$M$ Gorenstein weak injective?}
\end{question}

We call an $R$-module $M$ defined as  above  \emph{two-degree
Gorenstein weak injective}, and denote by $\mathcal{GWI}^2({R})$
the class of two-degree Gorenstein weak injective $R$-modules. It is
obvious that there is a containment $\mathcal{GWI}(R)\subseteq
\mathcal{GWI}^2({R})$.

In the following, we show that the answer to the above question is
affirmative over  our hypothesis that the class $\mathcal{GWI}(R)$ is closed under extensions.

\begin{theorem}\label{main}
$\mathcal{GWI}(R)=\mathcal{GWI}^2({R})$.
\end{theorem}

Before giving the proof of Theorem $\ref{main}$, we need the
following preliminaries.

\begin{definition}\label{strongly two-degree}
An $R$-module $M$ is called  \emph{strongly two-degree Gorenstein
weak injective} if there exists an exact sequence
$$
\xymatrix@C=0.5cm{
  \cdots \ar[r] & G \ar[r]^{d} & G \ar[r]^{d} & G \ar[r]^{d} & G \ar[r] & \cdots },
$$
where $G$ is Gorenstein weak injective, such that $M=\mbox{Coker}d$
and the functor $\operatorname{Hom}_R(N,-)$ leaves this sequence
exact whenever $N$ is a super finitely presented $R$-module with
$pd_R(N)<\infty$.
\end{definition}

We denote by $\mathcal{SGWI}^2({R})$ the class of strongly
two-degree Gorenstein weak injective $R$-modules. It is obvious that
there is a containment $
\mathcal{SGWI}^2({R})\subseteq\mathcal{GWI}^2({R})$. As a similar argument to the proof of  Proposition \ref{simplify}, we have

\begin{lemma}\label{stwochar}
Let $M$ be an $R$-module. Then the following are equivalent:
\begin{itemize}
  \item [(1)] $M$ is strongly two-degree Gorenstein weak injective;
  \item [(2)] There exists an exact sequence
$$
\xymatrix@C=0.5cm{
  \cdots \ar[r] & G \ar[r]^{d} & G \ar[r]^{d} & G \ar[r]^{d} & G \ar[r] & \cdots },
$$
where $G$ is Gorenstein weak injective;
  \item [(3)] There exists an exact sequence $0\rightarrow M\rightarrow
G\rightarrow M\rightarrow 0$, where $G$ is Gorenstein weak
injective.
\end{itemize}
\end{lemma}

\begin{proposition} Let $M$ be an $R$-module.
If $M$ is  two-degree Gorenstein weak injective, then it is a direct
summand of some strongly two-degree Gorenstein weak injective
$R$-module.
\end{proposition}

\begin{proof}
Since $M$ is  two-degree Gorenstein weak injective, there exists an
exact sequence of Gorenstein weak injective $R$-modules
$$
\mathbb{G}=\xymatrix@C=0.5cm{
  \cdots \ar[r] & G_1 \ar[r]^{d_1} & G_0 \ar[r]^{d_0} & G_{-1} \ar[r]^{d_{-1}} & G_{-2} \ar[r] & \cdots }
$$
where $G_{-i}=G^{i-1}$ for each $i\geq 1$, such that
$M=\mbox{Im}d_0$. For all $m\in \mathbb{Z}$, we
denote by $\Sigma ^m\mathbb{G}$ the exact sequence obtained
from $\mathbb{G}$ by increasing all indexes by $m$: $(\Sigma
^m\mathbb{G})_i=G_{i-m}$ and $d^{\Sigma
^m\mathbb{G}}_i=d_{i-m}$ for all $i\in \mathbb{Z}$. Then we
get the following exact sequence
$$
\bigoplus_{m\in \mathbb{Z}} (\Sigma
^m\mathbb{G})=\xymatrix@C=0.5cm{
  \cdots \ar[r] & \bigoplus_{i\in\mathbb{Z}} G_i \ar[r]^{d} & \bigoplus_{i\in\mathbb{Z}} G_i \ar[r]^{d} & \bigoplus_{i\in\mathbb{Z}} G_i \ar[r]^{d} & \cdots },
$$
where $d([g_i])=[d_i(g_i)]$ for any $g_i\in G_i$. Then $d^2=0$ and
$\bigoplus_{i\in\mathbb{Z}} G_i$ is Gorenstein weak injective  by
Remark $\ref{remark1}(3)$. Let $H=\mbox{Im}d$. Then $H\subseteq
\mbox{Ker}d$. Now we assume that $d([g_i])=[d_i(g_i)]=0$. It is
clear that $d_i(g_i)=0$, and so there exists $g_{i+1}\in G_{i+1}$
such that $d_{i+1}(g_{i+1})=g_i$. Thus $d([g_{i+1}])=[g_i]$, and
hence $\mbox{Ker}d\subseteq \mbox{Im}d=H$. Therefore, we may get an
exact sequence $0\rightarrow H\rightarrow G\rightarrow H\rightarrow
0$, where $G:=\bigoplus_{i\in\mathbb{Z}} G_i$.  Thus
$H$ is strongly two-degree Gorenstein weak injective by Lemma
$\ref{stwochar}$.

Finally, since $H\cong \oplus_{i\in\mathbb{Z}}\mbox{Im}d_i$, we get
that $M$ is a direct summand of $H$, as desired.
\end{proof}

Recall from \cite[1.1]{Ho} that a class $\mathcal{C}$ of $R$-modules
is called \emph{injectively resolving} if all injective $R$-modules
are contained in $\mathcal{C}$, and for any exact sequence
$0\rightarrow L\rightarrow M\rightarrow N\rightarrow 0$ with $L\in
\mathcal{C}$, the conditions $M\in \mathcal{C}$ and $N\in
\mathcal{C}$ are equivalent.

\begin{proposition}\label{injectively resolving}
The class $\mathcal{GWI}(R)$ is
injectively resolving.
\end{proposition}

\begin{proof}
Clearly, every injective $R$-module is Gorenstein weak injective.
Moreover, the class $\mathcal{GWI}(R)$  is closed under extensions, by our running hypothesis. So we will prove that for any exact sequence
$0\rightarrow L\rightarrow M\rightarrow N\rightarrow 0$, if  $L$ and
$M$ are Gorenstein weak injective, then so is $N$. Indeed, since $M$
is Gorenstein weak injective, there exists an exact sequence
$0\rightarrow G\rightarrow W\rightarrow M\rightarrow 0$ such that
$W$ is weak injective and $G$ is Gorenstein weak injective.

Consider the following pull-back diagram:
\begin{gather}\label{p}
\begin{split}
\xymatrix@=20pt{
&0\ar[d]&0\ar[d]&&\\
&G\ar@{=}[r]\ar[d]&G\ar[d]&&\\
0\ar[r]&N'\ar[r]\ar[d]&W\ar[r]\ar[d]&N\ar@{=}[d]\ar[r]&0\\
0\ar[r]&L\ar[r]\ar[d]&M\ar[r]\ar[d]&N\ar[r]&0\\
&0&0&& }
\end{split}
\end{gather}
Since $L$ and $G$ are Gorenstein weak injective, $N'$ is also
Gorenstein weak injective. Thus $N$ is Gorenstein weak
injective by Proposition $\ref{simplify}$.
\end{proof}

\begin{corollary}\label{direct summands}
The class $\mathcal{GWI}(R)$
is closed under direct summands.
\end{corollary}

\begin{proof}
It follows from \cite[Prop. 1.4]{Ho}, Remark $\ref{remark1}(3)$ and
Proposition $\ref{injectively resolving}$.
\end{proof}

Now we give the proof of our main theorem (Theorem $\ref{main}$).

\medskip

\emph{Proof of Theorem $\ref{main}$.}  Since
$\mathcal{GWI}({R})\subseteq \mathcal{GWI}^2({R})$, it suffices to
show that $\mathcal{GWI}^2({R})\subseteq \mathcal{GWI}({R})$. Since
every two-degree Gorenstein weak injective $R$-module is a direct
summand of some strongly two-degree Gorenstein weak injective
$R$-module, and the class of Gorenstein weak injective $R$-modules
is closed under direct summands by Corollary $\ref{direct
summands}$, so we only need to prove that every strongly two-degree
Gorenstein weak injective $R$-module is Gorenstein weak injective.

 Let $M$ be a strongly two-degree Gorenstein weak
injective $R$-module. Then, by Lemma $\ref{stwochar}$, there is an
exact sequence $0\rightarrow M\rightarrow G\rightarrow M\rightarrow
0$ with $G$ Gorenstein weak injective. Moreover, there is an exact
sequence $0\rightarrow G_1\rightarrow W\rightarrow G\rightarrow 0$
with $W$ weak injective and $G_1$ Gorenstein weak injective.

Consider the following pull-back diagram:
\begin{gather}\label{6}
\begin{split}
\xymatrix@=20pt{
&0\ar[d]&0\ar[d]&&\\
&G_1\ar@{=}[r]\ar[d]&G_1\ar[d]&&\\
0\ar[r]&N\ar[r]\ar[d]&W\ar[r]\ar[d]&M\ar@{=}[d]\ar[r]&0\\
0\ar[r]&M\ar[r]\ar[d]&G\ar[r]\ar[d]&M\ar[r]&0\\
&0&0&& }
\end{split}
\end{gather}
From the middle row in the diagram ($\ref{6}$), we obtain an exact
sequence $0\rightarrow N\rightarrow W\rightarrow M\rightarrow 0$
with $W$ weak injective. Thus, in order to show that $M$ is
Gorenstein weak injective, it suffices to prove that $N$ is
Gorenstein weak injective by  Proposition $\ref{simplify}$.

Consider the following pull-back diagram:
\begin{gather}\label{7}
\begin{split}
\xymatrix@=20pt{
&&0\ar[d]&0\ar[d]&\\
&&G_1\ar@{=}[r]\ar[d]&G_1\ar[d]&\\
0\ar[r]&M\ar[r]\ar@{=}[d]&G_2\ar[r]\ar[d]&N\ar[d]\ar[r]&0\\
0\ar[r]&M\ar[r]&G\ar[r]\ar[d]&M\ar[r]\ar[d]&0\\
&&0&0& }
\end{split}
\end{gather}
Since $G$ and $G_1$ are Gorenstein weak injective, $G_2$ is also
Gorenstein weak injective by the middle column in the diagram
($\ref{7}$). Hence there exists an exact sequence $0\rightarrow
G_3\rightarrow W_0\rightarrow G_2\rightarrow 0$ such that $W_0$ is
weak injective and $G_3$ is Gorenstein weak injective.

 Consider the following pull-back diagram:
\begin{gather}\label{8}
\begin{split}
\xymatrix@=20pt{
&0\ar[d]&0\ar[d]&&\\
&G_3\ar@{=}[r]\ar[d]&G_3\ar[d]&&\\
0\ar[r]&N_1\ar[r]\ar[d]&W_0\ar[r]\ar[d]&N\ar@{=}[d]\ar[r]&0\\
0\ar[r]&M\ar[r]\ar[d]&G_2\ar[r]\ar[d]&N\ar[r]&0\\
&0&0&& }
\end{split}
\end{gather}
From the middle row in the diagram ($\ref{8}$), we obtain an exact
sequence  $0\rightarrow N_1\rightarrow W_0\rightarrow N\rightarrow
0$ with $W_0$ weak injective. We repeat the argument by replacing
$N$ with $N_1$ to get $N_2$ and an exact sequence $0\rightarrow
N_2\rightarrow W_1\rightarrow N_1\rightarrow 0$ with $W_1$ weak
injective. Continuing this process, we may obtain an exact sequence
$\cdots \rightarrow W_1\rightarrow W_0\rightarrow N\rightarrow 0$,
where each $W_i$ is weak injective, which shows that $N$ is
Gorenstein weak injective by Proposition $\ref{simplify}$. We have
completed the proof.

\section*{Acknowledgement}
This work was partially
supported by NSFC (Grant Nos. 11571164, 11571341).
The authors thank the referee for the useful suggestions.

\end{document}